# RANDOM REWARDS, FRACTIONAL BROWNIAN LOCAL TIMES AND STABLE SELF-SIMILAR PROCESSES

By Serge Cohen and Gennady Samorodnitsky[1]

*Université Paul Sabatier and Cornell University*

We describe a new class of self-similar symmetric $\alpha$-stable processes with stationary increments arising as a large time scale limit in a situation where many users are earning random rewards or incurring random costs. The resulting models are different from the ones studied earlier both in their memory properties and smoothness of the sample paths.

**1. Introduction.** With the dramatic increase of importance of communication networks came the need to better understand their behavior at different scales. This requires a construction of stochastic models that can plausibly arise as the result of activities associated with such networks. Limiting stochastic processes often scale, and one hopes that such models can provide insight into the scaling of properties of the networks.

Perhaps the best known result of this type is the paper [21], where the limiting model turned out to be (depending on the relationship between the number of users and the time scale) either fractional Brownian motion or a Lévy $\alpha$-stable motion (this paper followed up and was an improvement of the earlier papers [38] and [36]). The fact that either a light-tailed but long-range dependent model or a heavy-tailed but short-range dependent model could appear has become an article of faith; see, for example, [7] for an application in a network context. Other heavy-tailed limiting models have appeared (see, e.g., [25]), but the limiting processes are not long-range dependent (more about this will appear in the sequel).

In this paper we exhibit a natural situation where the limiting model belongs to a new class of $\alpha$-stable models. It is a self-similar process with

Received August 2005; revised January 2006.
[1]Supported in part by NSF Grant DMS-03-03493 and NSA Grant MSPF-02G-183 at Cornell University.

stationary increments, and we will argue that the increments are long-range dependent. Let $(W_k, k \in \mathbf{Z})$ be a sequence of i.i.d. random symmetric variables satisfying

$$\overline{F}_W(x) := P(W_0 > x) \sim \sigma_W^\alpha x^{-\alpha} \tag{1.1}$$

as $x \to \infty$, where $0 < \alpha < 2$ and $\sigma_W > 0$. Further, let $(V_1, V_2, \ldots)$ be a sequence of i.i.d. mean zero and unit variance integer-valued random variables, independent of $(W_k, k \in \mathbf{Z})$, defining a random walk $S_n = V_1 + \cdots + V_n$ for $n \geq 1$. If one views $S_n$ as describing the "position" of the "state" of a user at time $n$, and $W_k$ the "reward" earned by, "cost" incurred by or "amount of work" produced by the user in state $k$, then the total reward earned by time $n$ is

$$R(n) = \sum_{j=1}^{n} W_{S_j}. \tag{1.2}$$

Assuming that there are many such users earning independent rewards or generating independent work, it turns out that a properly normalized sequence of rewards converges weakly to a limit, which we will call an FBM-1/2-local time fractional symmetric $\alpha$-stable motion. This is a particular case of a larger class of models, FBM-$H$-local time fractional symmetric $\alpha$-stable motion, where $0 < H < 1$ (these are self-similar with exponent of self-similarity $H' = 1 - H + H/\alpha$). We will represent this process as a stochastic integral with respect to an $\alpha$-stable random measure, with the integrand being the local time process of a fractional Brownian motion with exponent $H$ (hence the name of the model). The increments of this process are generated by a conservative null flow (see below for the details) and, hence, this process turns out to be different from all other classes of $\alpha$-stable self-similar processes with stationary increments that have been considered so far in the extensive literature on the subject.

Two remarks have to be made at this point. First of all, the only reason for assuming symmetry of $(W_k, k \in \mathbf{Z})$ is that dealing with symmetric $\alpha$-stable (S$\alpha$S) models leads to simpler expressions and unified exposition for all $0 < \alpha < 2$. Classes of nonsymmetric stable models parallel to those we are working with in this paper can be defined without difficulty, the case $\alpha = 1$ being the exception. Under suitable tail conditions, the random reward scheme with appropriate translation and scaling will converge to these stable processes. Second, our processes are related to a family of limiting models obtained in similar circumstances (but with a single user) in [15]. In their case, the limiting process is self-similar with stationary increments, but not stable.

This paper is organized as follows. In the next section we will summarize the required information on $\alpha$-stable processes and random measures, on self-similarity and on local times of fractional Brownian motions. Our process is formally introduced in Section 3. The properties of the increment process are discussed in Section 4. In Section 5 we study the smoothness of the sample paths of local time fractional stable motions through their Hölder continuity properties. It turns out that local time fractional stable motions can be naturally written as sums of *absolutely continuous* self-similar stable processes, and the decomposition goes through the chaos expansion of the local times of fractional Brownian motions. This is done in Section 6. In Section 7 we prove the aforementioned convergence of the random reward scheme to the FBM-1/2-local time fractional stable motion. We conclude with some comments and a discussion of possible extensions in Section 8.

**2. Preliminaries.** Throughout this paper, we will deal with S$\alpha$S processes given in the form

$$(2.1) \qquad X(t) = \int_E f(t,x) M(dx), \qquad t \in T,$$

where $T$ is a parameter space, $M$ is a symmetric $\alpha$-stable random measure on a measurable space $(E, \mathcal{E})$ with a $\sigma$-finite control measure $m$, and $f(t, \cdot) \in L^\alpha(m, \mathcal{E})$ for all $t \in T$. See Chapter 3 of [34] for information on $\alpha$-stable random measures and integrals with respect to these measures.

If the parameter space is countable (e.g., if $T = \mathbf{Z}$), and the process is stationary (under the usual left shift operator), then it has an integral representation as above, but the kernels $f_n = f(n, \cdot)$, $n \in \mathbf{Z}$, are of a special form. Specifically, one can choose

$$(2.2) \qquad f_n(x) = a_n(x) \left( \frac{dm \circ \phi^n}{dm}(x) \right)^{1/\alpha} f \circ \phi^n(x), \qquad x \in E,$$

for $n = 0, 1, 2, \ldots$, where $\phi: E \to E$ is a measurable nonsingular map (i.e., a one-to-one map with both $\phi$ and $\phi^{-1}$ measurable, mapping the control measure $m$ into an equivalent measure), where

$$a_n(x) = \prod_{j=0}^{n-1} u \circ \phi^j(x), \qquad x \in E,$$

for $n = 0, 1, 2, \ldots$, with $u: E \to \{-1, 1\}$ a measurable function and $f \in L^\alpha(m, \mathcal{E})$; see [29]. Many properties of the resulting stable process are closely connected with the ergodic-theoretic properties of the flow (the group of maps) $(\phi^n, n \in \mathbf{Z})$, an important classification of which is into dissipative, conservative null and positive flows; see [29, 32] and [33]. In particular, a key idea in the latter two papers is that it is possible to view stationary stable processes corresponding to dissipative flows as short memory processes, those corresponding to positive flows as infinite memory processes and those corresponding to conservative null flows as processes with a finite but long memory. Good general references on ergodic theory are [1] and [17].

A stochastic process $(Y(t), t \geq 0)$ is called *self-similar* with *exponent of self-similarity* $H$ if, for all $c > 0$, the processes $(Y(ct), t \geq 0)$ and $(c^H Y(t), t \geq 0)$ have the same finite-dimensional distributions. Most commonly studied are self-similar processes with stationary increments $[(Y(t+a) - Y(a), t \geq 0)$ has the same finite-dimensional distributions for all $a \geq 0]$. The common abbreviation for such a process is SSSI (self-similar stationary increments), or $H$-SSSI, if the exponent of self-similarity $H$ is to be emphasized.

For SSSI processes with a finite mean, the exponent of self-similarity is restricted to the range $0 < H < 1$ (apart from degenerate cases) and, in that range, there is a unique $H$-SSSI Gaussian process. It has zero mean and covariance function

$$\mathrm{Cov}(Y(s), Y(t)) = \frac{EY^2(1)}{2}[t^{2H} + s^{2H} - (t-s)^{2H}],$$

$0 \leq s \leq t$. This process is called the fractional Brownian motion (FBM).

In the $\alpha$-stable case, $0 < \alpha < 2$, the family of SSSI processes is much larger. The feasible range of pairs $(\alpha, H)$ is

(2.3) $$\begin{cases} 0 < H \leq 1/\alpha, & \text{if } 0 < \alpha \leq 1, \\ 0 < H < 1, & \text{if } 1 < \alpha < 2, \end{cases}$$

and, apart from the case $0 < \alpha < 1$ and $H = 1/\alpha$, a feasible pair $(\alpha, H)$ does not determine the law of an S$\alpha$S $H$-SSSI process.

REMARK 2.1. The class of S$\alpha$S SSSI processes constructed in this paper has exponent of self-similarity in the range

(2.4) $$\begin{cases} 1 < H < 1/\alpha, & \text{if } 0 < \alpha < 1, \\ H = 1, & \text{if } \alpha = 1, \\ 1/\alpha < H < 1, & \text{if } 1 < \alpha < 2. \end{cases}$$

It has been a long-standing challenge to describe classes of symmetric 1-stable SSSI processes with $H = 1$ other than linear combinations of independent symmetric 1-stable Lévy motion and the straight line process $Y(t) = tY(1), t \geq 0$. The model developed in this paper provides, in the particular case $\alpha = 1$, an entire family of such processes.

Two of the most well-known families of S$\alpha$S $H$-SSSI processes (with $0 < H < 1$) are obtained by taking two of the many possible integral representations of the fractional Brownian motion and modifying them appropriately (in particular, replacing the Brownian motion as an integrator by an S$\alpha$S Lévy motion). These are the linear fractional stable motion and the real harmonizable fractional stable motion. The linear fractional stable motion belongs to the class of self-similar stable mixed moving averages described by [26, 27]. Its increment process is generated by a dissipative flow. On the other hand, the increment process of the real harmonizable fractional stable motion is generated by a positive flow. We refer the reader to Chapter 7 of [34] and to [10] for more information on self-similar processes.

A fractional Brownian motion with any exponent of self-similarity $0 < H < 1$ has a local time process $(l(x,t), x \in \mathbb{R}, t \geq 0)$ that is jointly continuous in $x$ and $t$ [3]. The self-similarity property of the fractional Brownian motion immediately implies the scaling property of the local time process: for any $c > 0$,

$$(2.5) \qquad (l(c^H x, ct), x \in \mathbb{R}, t \geq 0) \stackrel{d}{=} (c^{1-H} l(x,t), x \in \mathbb{R}, t \geq 0),$$

a somewhat more convenient form of which is

$$(2.6) \qquad \left(\frac{1}{c} l(x, ct), x \in \mathbb{R}, t \geq 0\right) \stackrel{d}{=} \left(\frac{1}{c^H} l\left(\frac{x}{c^H}, t\right), x \in \mathbb{R}, t \geq 0\right).$$

It is a simple consequence of (2.5) and of Theorem 6, page 275 of [14] that, on a set of probability 1,

$$(2.7) \qquad \lim_{t \to \infty} l(x,t) = \infty \qquad \text{for all } x \in \mathbb{R}.$$

Similarly, the stationarity of increments property of the fractional Brownian motion implies a type of stationarity of the increments of the local time, which can be formulated as follows. Let $(\Omega, \mathcal{F}, \mathbf{P})$ be the probability space on which the fractional Brownian motion and its local time process live. Then, abusing somewhat the term "law" by applying it to an infinite induced measure,

$$(2.8) \qquad \begin{aligned} &\text{the law of } (l(x, t+h)(\omega) - l(x, h)(\omega), t \geq 0) \\ &\text{under } \mathbf{P} \times \text{Leb does not depend on } h \geq 0. \end{aligned}$$

A modification of the proof of Theorem 1.2 that leads to Corollary 1.1 in [39] gives us that

$$(2.9) \qquad K := \sup_{\substack{x \in \mathbb{R} \\ 0 \leq s < t \leq 1/2}} \frac{l(x,t) - l(x,s)}{(t-s)^{1-H}(\log 1/(t-s))^H} < \infty \qquad \text{a.s.}$$

and has finite moments of all orders. (Note that using instead the estimates in [6] gives a slightly worse power of the logarithm: $H+1$ instead of $H$.) In particular, $l(x,t)$ has moments of all orders finite and uniformly bounded in all real $x$ and all $t$ in a compact set.

**3. FBM-$H$-local time fractional stable motions.** We now introduce our class of models. Let $(\Omega', \mathcal{F}', \mathbf{P}')$ be a probability space supporting a fractional Brownian motion $(B_H(t), t \geq 0)$ with exponent of self-similarity $H$, and let $l = l(x,t) = l(x,t)(\omega')$ be its jointly continuous local time process. Let $M$ be an S$\alpha$S random measure on the space $\Omega' \times \mathbb{R}$ with control measure $\mathbf{P}' \times \text{Leb}$, where Leb is the Lebesgue measure on $\mathbb{R}$. The random measure itself lives on some other probability space $(\Omega, \mathcal{F}, \mathbf{P})$. We define

$$(3.1) \qquad Y(t) = \int_{\Omega'} \int_{\mathbb{R}} l(x,t)(\omega') M(d\omega', dx), \qquad t \geq 0.$$

Our first result below shows that $(Y(t), t \geq 0)$ is a well-defined S$\alpha$S process which is self-similar and has stationary increments. We call this process *FBM-$H$-local time fractional symmetric $\alpha$-stable motion.*

THEOREM 3.1. *The process $(Y(t), t \geq 0)$ in (3.1) is a well-defined S$\alpha$S process. It has stationary increments and is self-similar, with exponent*

$$(3.2) \qquad H' = 1 - H + H/\alpha = 1 + H\left(\frac{1}{\alpha} - 1\right).$$

PROOF. To show that $Y$ is properly defined we need to check that

$$\int_{\Omega'} \int_{\mathbb{R}} l^{\alpha}(x,t)(\omega') \mathbf{P}'(d\omega') \, dx = \mathbf{E}' \int_{\mathbb{R}} l^{\alpha}(x,t) \, dx < \infty.$$

We have

$$(3.3) \qquad \begin{aligned} \mathbf{E}' \int_{\mathbb{R}} l^{\alpha}(x,t) \, dx &= \int_{\mathbb{R}} \mathbf{E}' \left[ l^{\alpha}(x,t) \mathbf{1}\left( \sup_{0 \leq s \leq t} |B_H(s)| \geq |x| \right) \right] dx \\ &\leq \int_{\mathbb{R}} (\mathbf{E}' l^2(x,t))^{\alpha/2} \left( \mathbf{P}'\left( \sup_{0 \leq s \leq t} |B_H(s)| \geq |x| \right) \right)^{1/q} dx \end{aligned}$$

with $q = 1 - \alpha/2$. Since the moments of the local time are uniformly bounded and

$$(3.4) \qquad \int_{\mathbb{R}} \left( \mathbf{P}'\left( \sup_{0 \leq s \leq t} |B_H(s)| \geq |x| \right) \right)^{1/q} dx < \infty$$

as the supremum of a bounded Gaussian process has Gaussian-like tails, we conclude that the left-hand side of (3.3) is finite and, hence, $(Y(t), t \geq 0)$ in (3.1) is a well-defined S$\alpha$S process.

Notice that for any $c > 0, k \geq 1, \theta_1, \ldots, \theta_k \in \mathbb{R}$ and $t_1, \ldots, t_k \geq 0$ we have, using (2.6),

$$\mathbf{E} \exp\left( i \sum_{j=1}^{k} \theta_j Y(ct_j) \right) = \exp\left( -\int_{\mathbb{R}} \mathbf{E}' \left| \sum_{j=1}^{k} \theta_j l(x, ct_j) \right|^{\alpha} dx \right)$$

$$= \exp\left(-\int_{\mathbb{R}} \mathbf{E}' \left|\sum_{j=1}^{k} \theta_j c^{1-H} l\left(\frac{x}{c^H}, t_j\right)\right|^\alpha dx\right)$$

$$= \exp\left(-c^{\alpha(1-H)} \mathbf{E}' \int_{\mathbb{R}} \left|\sum_{j=1}^{k} \theta_j l\left(\frac{x}{c^H}, t_j\right)\right|^\alpha dx\right)$$

$$= \exp\left(-c^{\alpha(1-H)+H} \mathbf{E}' \int_{\mathbb{R}} \left|\sum_{j=1}^{k} \theta_j l(y, t_j)\right|^\alpha dy\right)$$

$$= \mathbf{E} \exp\left(i \sum_{j=1}^{k} \theta_j c^{1-H+H/\alpha} Y(t_j)\right).$$

Therefore, $(Y(t), t \geq 0)$ is $H'$-self-similar, with $H'$ given by (3.2).

Furthermore, for any $h \geq 0, k \geq 1, \theta_1, \ldots, \theta_k \in \mathbb{R}$ and $t_1, \ldots, t_k \geq 0$ we have, by (2.8)

(3.5)
$$\mathbf{E} \exp\left(i \sum_{j=1}^{k} \theta_j (Y(t_j + h) - Y(h))\right)$$
$$= \exp\left(-\int_{\mathbb{R}} \mathbf{E}' \left|\sum_{j=1}^{k} \theta_j (l(x, t_j + h) - l(x, h))\right|^\alpha dx\right)$$
$$= \exp\left(-\int_{\mathbb{R}} \mathbf{E}' \left|\sum_{j=1}^{k} \theta_j l(x, t_j)\right|^\alpha dx\right)$$
$$= \mathbf{E} \exp\left(i \sum_{j=1}^{k} \theta_j Y(t_j)\right).$$

Therefore, $(Y(t), t \geq 0)$ has stationary increments. $\square$

REMARK 3.2. Observe that:

1. For $0 < \alpha < 1$ we obtain a family of $H'$-SSSI S$\alpha$S processes with $H' \in (1, 1/\alpha)$.
2. For $1 < \alpha < 2$ we obtain a family of $H'$-SSSI S$\alpha$S processes with $H' \in (1/\alpha, 1)$.
3. For $\alpha = 1$ we obtain a family of 1-SSSI S$\alpha$S processes.

Notice that, for $\alpha \neq 1$, different choices of the fractional Brownian motion exponent of self-similarity $H$ lead to a different exponent of self-similarity $H'$ of the S$\alpha$S process $(Y(t), t \geq 0)$ and, hence, to a different process. On the other hand, for $\alpha = 1$ the exponent of self-similarity $H'$ is *independent* of $H$. Nonetheless, the processes $(Y(t), t \geq 0)$ corresponding to different $H$ are different in this case as well, as will be seen in the sequel.

**4. The increment process.** An object of interest for an SSSI process is its increment process. It is a stationary process, and its memory properties are often of interest. For example, the increment process of a fractional Brownian motion, the so-called fractional Gaussian noise, is a standard long memory (if $H > 1/2$) model that was used by Mandelbrot (see, e.g., [19, 20]) to explain the famous Hurst phenomenon. Similarly, the increments of the linear fractional stable motion are called linear fractional stable noise, and those of the real harmonizable fractional stable motion are called (real) harmonizable fractional stable noise. It is often believed that the properties of the fractional noises are largely determined by the exponent of self-similarity of the original process. One of the goals of this section (which studies the increment process of the FBM-$H$-local time fractional $\alpha$-stable motion) is to shed some light on this question.

Let, therefore, $(Y(t), t \geq 0)$ be an FBM-$H$-local time fractional S$\alpha$S motion and consider its increment process

$$(4.1) \qquad Z_n = Y(n+1) - Y(n), \qquad n = 0, 1, \ldots,$$

which will be called *FBM-$H$-local time fractional S$\alpha$S noise.*

A very important property of the FBM-$H$-local time fractional S$\alpha$S noise is given in the following result:

THEOREM 4.1. *The FBM-$H$-local time fractional S$\alpha$S noise is generated by a conservative null flow.*

PROOF. Note that the FBM-$H$-local time fractional S$\alpha$S noise has an integral representation

$$(4.2) \quad Z_n = \int_{\Omega'} \int_{\mathbb{R}} (l(x, n+1)(\omega') - l(x, n)(\omega')) M(d\omega', dx), \qquad n \geq 0.$$

Let $C$ be the space of continuous functions from $\mathbb{R}$ to $\mathbb{R}$ and $\mathbf{P}'_1$ a probability measure on $C$, under which the coordinate map is the fractional Brownian motion with exponent of self-similarity $H$. Let $m$ be a $\sigma$-finite measure on $C$ defined by $m = (\mathbf{P}'_1 \times \text{Leb}) \circ T^{-1}$, where $T: C \times \mathbb{R} \to C$ is given by $T(\omega', x) = \omega' - x$, $\omega' \in C, x \in \mathbb{R}$. Let $L: C \to \mathbb{R}$ be a measurable function that associates to a function $\omega' \in C$ its local time at 0 in the interval $(0, 1]$ if $\omega'$ has continuous local time. An alternative representation of the process in (4.2) is then

$$(4.3) \qquad Z_n = \int_C L \circ \phi^n(\omega') M_1(d\omega'), \qquad n \geq 0,$$

where $M_1$ is an S$\alpha$S random measure on $C$ with control measure $m$, and $\phi: C \to C$ is given by $\phi(\omega') = \omega'(\cdot + 1)$. The stationarity of the increments of the fractional Brownian motion implies that the map $\phi$ preserves the measure $m$. Note that (4.3) is a representation of type (2.2) (with both $a_n \equiv 1$ and the Radon–Nikodym derivative equal to 1). A conclusion is that the flow $(\phi^n)$ and the underlying measure space on which $(\phi^n)$ acts are the same, *independently of the value of $\alpha$*. Therefore, it is sufficient to prove the theorem in the case $\alpha = 1$, which we will assume until the end of the proof.

We continue working with the representation (4.2). Note that, by (2.7),

$$\text{(4.4)} \quad \sum_{n=0}^{m} [l(x, n+1)(\omega') - l(x, n)(\omega')]$$
$$= l(x, m+1)(\omega') \to \infty \quad \text{as } m \to \infty$$

outside a subset of $\Omega' \times \mathbb{R}$ of measure 0. By Corollary 4.2 of [29], this implies that the FBM-$H$-local time fractional S$\alpha$S noise is generated by a conservative flow. It also, evidently, shows that the kernel in the representation (4.2) has a full support.

In order to prove that the FBM-$H$-local time fractional S$\alpha$S noise is generated by a null flow, we will apply Corollary 2.2 of [33] to the obvious two-sided extension of the process to $(Z_n, n \in \mathbf{Z})$. For reasons of symmetry, it is enough to exhibit a nonincreasing nonnegative sequence $w_n$ such that

$$\text{(4.5)} \quad \sum_{n=0}^{\infty} w_n = \infty$$

and

$$\text{(4.6)} \quad \sum_{n=0}^{\infty} w_n [l(x, n+1)(\omega') - l(x, n)(\omega')] < \infty$$

for $\mathbf{P}' \times \text{Leb}$-almost every $(\omega', x)$.

Let $w_n = (1+n)^{-\theta}$ with some $1 - H < \theta \leq 1$. Since $\theta \leq 1$, the condition (4.5) is satisfied. To check (4.6), it is clearly enough to find a strictly positive measurable function $g$ such that

$$\text{(4.7)} \quad \mathbf{E}' \int_{\mathbb{R}} g(x) \sum_{n=0}^{\infty} w_n [l(x, n+1)(\omega') - l(x, n)(\omega')] \, dx < \infty.$$

Note that

$$\mathbf{E}' \int_{\mathbb{R}} g(x) \sum_{n=0}^{\infty} w_n [l(x, n+1)(\omega') - l(x, n)(\omega')] \, dx = \sum_{n=0}^{\infty} w_n \int_{n}^{n+1} \mathbf{E}' g(B_H(t)) \, dt.$$

Choose $g(x) = \exp(-x^2/2)$ so that for all $t \geq 0$

$$\mathbf{E}' g(B_H(t)) = \frac{1}{(1 + t^{2H}\sigma^2)^{1/2}},$$

where $\sigma^2 = \operatorname{Var} B_H(1)$. The left-hand side of (4.7) is then

$$(4.8) \qquad \sum_{n=0}^{\infty} w_n \int_n^{n+1} \frac{dt}{(1+t^{2H}\sigma^2)^{1/2}} \leq \sum_{n=0}^{\infty} w_n \frac{1}{(1+n^{2H}\sigma^2)^{1/2}} < \infty$$

by the choice of $\theta$. Hence, (4.6) follows. $\square$

REMARK 4.2. It follows from Theorem 4.1 that (for $1 < \alpha < 2$) the FBM-$H$-local time fractional S$\alpha$S motion is different from the linear fractional S$\alpha$S motion (or, more generally, from the self-similar mixed average processes of [26]) since the increments of the latter are generated by dissipative flows, and it is also different from the real harmonizable fractional S$\alpha$S motion whose increments are generated by positive flows.

In particular, the FBM-$H$-local time fractional S$\alpha$S noise can be viewed as a long memory process; its memory is longer than that of the linear fractional S$\alpha$S noise, but shorter than that of the harmonizable fractional S$\alpha$S noise. Implications of this will be seen, in particular, when we discuss smoothness of the sample paths in the next section. This is a reminder that very little is determined merely by the exponent of self-similarity for $\alpha$-stable SSSI processes.

REMARK 4.3. It follows immediately from Theorem 4.1 and Theorem 3.1 in [33] that the FBM-$H$-local time fractional S$\alpha$S noise is (unlike the harmonizable fractional S$\alpha$S noise) ergodic. It is easy to show that it is also a mixing process. Indeed, it suffices to show that, for any $0 < a < b$ and $\epsilon > 0$,

$$(4.9) \quad \begin{aligned} \lim_{n \to \infty} (\mathbf{P}' \times \mathrm{Leb}) \{ (\omega', x) : a \leq l(x,1)(\omega') \leq b, \\ l(x,n+1)(\omega') - l(x,n)(\omega') > \epsilon \} = 0; \end{aligned}$$

see, for example, [11] or [30]. Since the left-hand side of (3.3) is finite, we see that

$$(\mathbf{P}' \times \mathrm{Leb})\{(\omega', x) : a \leq l(x,1)(\omega') \leq b\} < \infty,$$

and so, given $\delta > 0$, for $K$ large enough,

$$(\mathbf{P}' \times \mathrm{Leb})\left\{ (\omega', x) : a \leq l(x,1)(\omega') \leq b, \sup_{0 \leq t \leq 1} |B_H(t)| > K \right\} \leq \delta.$$

For such $K$,

$$\begin{aligned} (\mathbf{P}' &\times \mathrm{Leb})\{(\omega', x) : a \leq l(x,1)(\omega') \leq b, l(x,n+1)(\omega') - l(x,n)(\omega') > \epsilon\} \\ &\leq \delta + (\mathbf{P}' \times \mathrm{Leb})\{(\omega', x) : |x| \leq K, l(x,n+1)(\omega') - l(x,n)(\omega') > \epsilon\} \\ &\leq \delta + 2K\mathbf{P}'(B_H(t) \in [-K,K] \text{ for some } n < t \leq n+1). \end{aligned}$$

Since the last probability clearly goes to zero as $n \to \infty$, we conclude that

$$\limsup_{n \to \infty} \mathbf{P}' \times \text{Leb}\{(\omega', x) : a \leq l(x,1)(\omega') \leq b,$$
$$l(x, n+1)(\omega') - l(x,n)(\omega') > \epsilon\} \leq \delta,$$

which proves (4.9), since $\delta > 0$ is arbitrary.

We close this section by addressing the point mentioned in Remark 3.2. Since, in the case $\alpha = 1$, the exponent of self-similarity of an FBM-$H$-local time fractional motion does not depend on $H$, one may suspect that $H$ does not change the law of the process itself (up to, perhaps, a multiplicative constant). The following result shows that this is not the case, and so the parameter $H$ parameterizes an entire family of different 1-stable SSSI processes that does not have either a Lévy 1-stable motion or a straight line process as a component (indeed, the former would have introduced a dissipative component to the flow generating the increment process, while the latter would have introduced a positive component to that flow):

PROPOSITION 4.4. *Let $\alpha = 1$ and $0 < H_1, H_2 < 1$ with $H_1 \neq H_2$. Then, there is no constant $C$ such that*

$$(Y_{H_1}(t), t \geq 0) \stackrel{\mathrm{d}}{=} (CY_{H_2}(t), t \geq 0),$$

*where $(Y_{H_i}(t), t \geq 0)$ is an FBM-$H_i$-local time fractional motion with $\alpha = 1$, $i = 1, 2$.*

PROOF. Assume that $H_1 < H_2$. If some $C$ with the above property existed, then we could use the fact that the kernel in the representation (4.2) has full support and Theorem 1.1 of [29] to connect the kernels with different $H$. Specifically, there would exist measurable maps

$$A : \Omega' \times \mathbb{R} \mapsto \mathbb{R} \setminus \{0\},$$
$$\Phi_1 : \Omega' \times \mathbb{R} \mapsto \mathbb{R},$$
$$\Phi_2 : \Omega' \times \mathbb{R} \mapsto \Omega',$$

such that

$$l_{H_1}(x, n+1)(\omega') - l_{H_1}(x, n)(\omega')$$
(4.10) $$= A(\omega', x)(l_{H_2}(\Phi_1(\omega', x), n+1)(\Phi_2(\omega', x))$$
$$- l_{H_2}(\Phi_1(\omega', x), n)(\Phi_2(\omega', x))), \qquad n \in \mathbf{N},$$

for $\mathbf{P}'_1 \times \text{Leb-almost every } \omega' \in \Omega', x \in \mathbb{R}$, where we have added subscripts to the local times with the obvious meaning, and $\mathbf{P}'_i$ is the probability measure on $\Omega'$ corresponding to the fractional Brownian motion with exponent $H_i$. Adding up, we obtain

$$(4.11) \quad l_{H_1}(x,n)(\omega') = A(\omega',x)l_{H_2}(\Phi_1(\omega',x),n)(\Phi_2(\omega',x)), \qquad n \in \mathbf{N},$$

for $\mathbf{P}'_1 \times \text{Leb}$-almost every $\omega' \in \Omega', x \in \mathbb{R}$.

By (2.6), Markov inequality and boundedness of the moments of the local time, there is a finite $K$ such that, for every $x \in \mathbb{R}, t > 0$ and $\varepsilon, \delta > 0$,

$$(4.12) \quad \mathbf{P}'_2(l_{H_2}(x,t) > \varepsilon t^{1-H_2+\delta}) \leq K\varepsilon^{-2}t^{-2\delta},$$

and so by Borel–Cantelli lemma

$$(4.13) \quad \mathbf{P}'_2(l_{H_2}(x,2^m) > \varepsilon 2^{m(1-H_2+\delta)} \text{ infinitely often in } m) = 0$$

for every $x \in \mathbb{R}$. By Fubini's theorem,

$$(4.14) \quad (\mathbf{P}'_2 \times \text{Leb})(G^c) = 0,$$

where

$$(4.15) \quad G = \left\{ (\omega',x), \lim_{m \to \infty} \frac{l_{H_2}(x,2^m)(\omega')}{2^{m(1-H_2+\delta)}} = 0 \right\}.$$

Therefore, in the definition (3.1) of the process $(Y(t), t \geq 0)$ for $H = H_2$, we can restrict the integral from $\Omega' \times \mathbb{R}$ to $G$ only and then, in (4.10) and (4.11), we will have

$$(4.16) \quad (\Phi_1(\omega',x), \Phi_2(\omega',x)) \in G$$

for all $\omega \in \Omega', x \in \mathbb{R}$. This means that for $\mathbf{P}'_1 \times \text{Leb}$-almost every $\omega \in \Omega', x \in \mathbb{R}$, we have

$$(4.17) \quad \lim_{m \to \infty} \frac{l_{H_1}(x,2^m)(\omega)}{2^{m(1-H_2+\delta)}} = 0.$$

Therefore, there is $x \in \mathbb{R}$, such that (4.17) holds $\mathbf{P}'_1$-a.s.

However, by (2.6),

$$\mathbf{P}'_1(l_{H_1}(x,2^m) > 2^{m(1-H_2+\delta)}) = \mathbf{P}'_1\left(l_{H_1}\left(\frac{x}{2^{mH_1}},1\right) > 2^{m(H_1-H_2+\delta)}\right).$$

If $\delta < H_2 - H_1$ this then gives us

$$\liminf_{m \to \infty} \mathbf{P}'_1(l_{H_1}(x,2^m) > 2^{m(1-H_2+\delta)}) \geq \mathbf{P}'_1(l_{H_1}(0,1) > 0) = 1 > 0,$$

contradicting (4.17). Therefore (4.10) is impossible and the proposition is proved. $\square$

**5. Hölder continuity.** The fact that the local times of the fractional Brownian motion are continuous and monotone in the time variable already implies that a FBM-$H$-local time fractional symmetric $\alpha$-stable motion with $0 < \alpha < 1$ is sample continuous (see, e.g., Theorem 10.4.2 of [34]) and the same is true for $1 < \alpha < 2$ by the mere fact that $H' > 1/\alpha$ (see Theorem 12.4.1 of [34]). Our goal in this section is to prove Hölder continuity of an FBM-$H$-local time fractional S$\alpha$S motion for all $0 < \alpha < 2$.

THEOREM 5.1. *Let $(Y(t), t \geq 0)$ be an FBM-$H$-local time fractional S$\alpha$S motion, $0 < \alpha < 2$. Then, it has a version with continuous sample paths satisfying*

$$(5.1) \qquad \sup_{0 \leq s < t \leq 1/2} \frac{|Y(t) - Y(s)|}{(t-s)^{1-H}(\log 1/(t-s))^{H+1/2}} < \infty \qquad a.s.$$

REMARK 5.2. It is instructive to express the Hölder continuity statement in (5.1) in terms of the exponent of self-similarity $H'$ of the FBM-$H$-local time fractional S$\alpha$S motion and $\alpha$, which can be done for $\alpha \neq 1$. For such $\alpha$, (5.1) means that an FBM-$H$-local time fractional S$\alpha$S motion is $d$-Hölder continuous with any

$$(5.2) \qquad d < \frac{H' - 1/\alpha}{1 - 1/\alpha}.$$

Let, for example, $1 < \alpha < 2$. Recall that a linear fractional S$\alpha$S motion with exponent of self-similarity $H' > 1/\alpha$ is $d$-Hölder continuous with any $d < H' - 1/\alpha$ [35], while a harmonizable fractional S$\alpha$S motion is $d$-Hölder continuous with any $d < H'$ [16]. In particular, an FBM-$H$-local time fractional S$\alpha$S motion has smoother sample paths than a linear fractional S$\alpha$S motion with the same exponent of self-similarity, and less smooth sample paths than a harmonizable fractional S$\alpha$S motion with the same exponent of self-similarity. This is not surprising if one recalls that the increments of an FBM-$H$-local time fractional S$\alpha$S motion have "stronger dependence" than those of a linear fractional S$\alpha$S motion, but not as strong as those of a harmonizable fractional S$\alpha$S motion.

Of course, since Theorem 5.1 only provides a lower bound on how smooth the sample functions are, the above discussion should be taken with "a grain of salt." We conjecture, however, that the upper bound on the Hölder exponent of the FBM-$H$-local time fractional S$\alpha$S motion cannot be improved. In the case $H = 1/2$, this is shown in Remark 5.3 below.

PROOF OF THEOREM 5.1. We will use a series representation of the stochastic integral (3.1) defining the FBM-$H$-local time fractional S$\alpha$S mo-tion. In distribution,

$$(5.3) \quad Y(t) = C_\alpha \sum_{j=1}^\infty G_j \Gamma_j^{-1/\alpha} e^{X_j^2/2\alpha} l_j(X_j, t), \qquad t \geq 0,$$

where $C_\alpha$ is a finite positive constant that depends only on $\alpha$, where $(G_j)$, $(\Gamma_j)$, $(X_j)$ and $(l_j)$ are four independent sequences such that $(G_j)$ and $(X_j)$ are i.i.d. standard normal random variables, $(\Gamma_j)$ are the arrival times of a unit rate Poisson process on $(0, \infty)$ and $(l_j)$ are i.i.d. copies of the local time process of a fractional Brownian motion; see Section 3.10 in [34].

Assume that the sequence $(G_j)$ is defined on some probability space $(\Omega_1, \mathcal{F}_1, \mathbf{P}_1)$, while the rest of the random variables on the right-hand side of (5.3) are defined on some other probability space $(\Omega_2, \mathcal{F}_2, \mathbf{P}_2)$, so that the FBM-$H$-local time fractional S$\alpha$S motion in the left-hand side of (5.3) is defined on the product of these two spaces. Let

$$K_j = \sup_{\substack{x \in \mathbb{R} \\ 0 \leq s < t \leq 1/2}} \frac{l_j(x,t) - l_j(x,s)}{(t-s)^{1-H}(\log 1/(t-s))^H}, \qquad j = 1, 2, \ldots,$$

and notice that, for a fixed $\omega_2 \in \Omega_2$, the process in (5.3) is centered Gaussian with the incremental variance

$$E_1(Y(t) - Y(s))^2 = C_\alpha^2 \sum_{j=1}^\infty \Gamma_j^{-2/\alpha} e^{X_j^2/\alpha} (l_j(X_j,t) - l_j(X_j,s))^2$$

$$(5.4) \qquad \leq \left( C_\alpha^2 \sum_{j=1}^\infty \Gamma_j^{-2/\alpha} e^{X_j^2/\alpha} K_j^2 \right) (t-s)^{2(1-H)} \left( \log \frac{1}{t-s} \right)^{2H}$$

$$:= M(\omega_2)(t-s)^{2(1-H)} \left( \log \frac{1}{t-s} \right)^{2H}$$

for all $0 \leq s < t \leq 1/2$, where $M$ is a $\mathbf{P}_2$-a.s. finite random variable on $(\Omega_2, \mathcal{F}_2, \mathbf{P}_2)$ (the latter statement follows from the fact that $E_2 K_j^\alpha < \infty$).

Applying now classical results on moduli of continuity of Gaussian processes (see, e.g., Theorem 2.1 of [8]) we obtain that, for $\mathbf{P}_2$-almost every $\omega_2 \in \Omega_2$,

$$\sup_{\substack{0 \leq s < t \leq 1/2 \\ s,t \text{ rational}}} \frac{|Y(t) - Y(s)|}{(t-s)^{1-H}(\log 1/(t-s))^{H+1/2}} < \infty, \qquad \mathbf{P}_1\text{-a.s.}$$

By Fubini's theorem,

$$\sup_{\substack{0 \leq s < t \leq 1/2 \\ s,t \text{ rational}}} \frac{|Y(t) - Y(s)|}{(t-s)^{1-H}(\log 1/(t-s))^{H+1/2}} < \infty, \qquad \mathbf{P}_1 \times \mathbf{P}_2\text{-a.s.},$$

which is equivalent to the statement of the theorem. □

REMARK 5.3. It is easy to show that, at least for $H = 1/2$, the result of Theorem 5.1 is "almost" sharp in the sense that there does not exist a function $g:(0,1/2) \to (0,\infty)$ with

$$\lim_{t \to 0} \frac{g(t)}{t^{1/2}(\log 1/t)^{1/2}} = 0$$

and, with positive probability,

(5.5) $$\sup_{0 \leq s < t \leq 1/2} \frac{|Y(t) - Y(s)|}{g(t-s)} < \infty.$$

Indeed, assume that such a function, in fact, exists. By the zero–one law for stable processes; see Section 9.5 in [34]. (5.5) would then hold with probability 1. It follows (e.g., by [28]), that we must have

(5.6) $$\sup_{0 \leq s < t \leq 1/2} \frac{|l(x,t) - l(x,s)|}{g(t-s)} < \infty \quad \text{a.s.}$$

for almost every $x \in \mathbb{R}$. Choose $x$ for which (5.6) holds and note that by the strong Markov property and [12],

$$\mathbf{P}'\left(\sup_{0 \leq s < t \leq 1/2} \frac{|l(x,t) - l(x,s)|}{g(t-s)} = \infty\right)$$
$$\geq \mathbf{P}'\left(\inf\{u \geq 0 : B_{1/2}(u) = x\} \leq \frac{1}{4}\right) \mathbf{P}'\left(\sup_{0 \leq s < t \leq 1/4} \frac{|l(0,t) - l(0,s)|}{g(t-s)} = \infty\right)$$
$$= \mathbf{P}'\left(\inf\{u \geq 0 : B_{1/2}(u) = x\} \leq \frac{1}{4}\right) > 0,$$

contradicting the necessity of (5.5) to hold with probability 1.

**6. Expansion into absolutely continuous terms.** The sample paths of (measurable) $H$-SSSI processes are almost never absolutely continuous with respect to the Lebesgue measure, the only exception being the case $H = 1$ with the process being the straight line process $Y(t) = tY(1)$ a.s. for all $t$ (see Theorem 3.3 of [37]). Nonetheless, there is a school of thought viewing nature as "producing smooth objects," with the others being more of a mathematical abstraction. In particular, smooth modifications of various mathematical models are of interest; see, for example, the "physical fractional Brownian motion" of [13]. In this section we use the chaos expansion of the local times of fractional Brownian motions due to [9] to construct an expansion of the FBM-$H$-local time fractional S$\alpha$S motion into a series of absolutely continuous S$\alpha$S self-similar processes, all with the same exponent of self-similarity as the original process. We first introduce the required notation.

For $\sigma > 0$, let $p_{\sigma^2}$ denote the density of a zero mean normal random variable with variance $\sigma^2$. $H_n$ is the $n$th Hermite polynomial

$$H_n(x) = \frac{(-1)^n}{n!} \exp\left(\frac{x^2}{2}\right) \frac{d^n}{dx^n}\left(\exp\left(-\frac{x^2}{2}\right)\right), \qquad x \in \mathbb{R},$$

with $H_0(x) \equiv 1$. Let $(\Omega', \mathcal{F}', \mathbf{P}')$ be a probability space supporting a Brownian motion $(W(s), s \in \mathbb{R})$, and let $I_n$ be the $n$th Wiener–Itô integral with respect to this Brownian motion. We refer the reader to [23] for information on these notions. Finally, let $K_H$ be the kernel defined by

$$K_H(t,s) = (t-s)^{H-1/2} - (H-1/2)\int_s^t (r-s)^{H-3/2}\left(1 - \left(\frac{s}{r}\right)^{-(H-1/2)}\right) dr$$

for $0 < s < t$ and equal to zero for other values of $s,t$. Note that, in distribution,

$$(6.1) \qquad B_H(t) = \left(\frac{\operatorname{Var} B_H(1)}{C_H}\right)^{1/2} \int_0^t K_H(t,s) W(ds), \qquad t \geq 0,$$

where $C_H$ is a finite positive constant depending only on $H$ (see, e.g., [2]).

THEOREM 6.1. *Let $(Y(t), t \geq 0)$ be an FBM-$H$-local time fractional $S\alpha S$ motion. In distribution,*

$$(6.2)\ Y(t) = \sum_{n=0}^{\infty} W_n(t) := \sum_{n=0}^{\infty} \int_{\Omega'} \int_{\mathbb{R}} h_n(x,t)(\omega') M(d\omega', dx), \qquad t \geq 0,$$

*where, for $n = 0, 1, \ldots,$*

$$(6.3) \quad \begin{aligned} h_n(x,t) &= h_n(x,t)(\omega') \\ &= \frac{C_H^{-n/2}}{\sigma} \int_0^t \frac{p_{s^{2H}}(x/\sigma)}{s^{nH}} H_n\left(\frac{x/\sigma}{s^H}\right) I_n(K_H(s,\cdot)^{\otimes n}) ds, \end{aligned}$$

*with $\sigma^2 = \operatorname{Var} B_H(1)$, and $M$ an $S\alpha S$ random measure on the space $\Omega' \times \mathbb{R}$ with control measure $\mathbf{P}' \times \operatorname{Leb}$. Each process $(W_n(t), t \geq 0)$ is a self-similar $S\alpha S$ process with exponent of self-similarity $H' = 1 - H + H/\alpha$ and has a modification with absolutely continuous sample paths. Moreover, the series in* (6.2) *converges in probability.*

PROOF. We first check that each $(W_n(t), t \geq 0)$ is a well-defined $S\alpha S$ process. Note that

$$(\sigma C_H^{n/2})^\alpha \int_{\Omega'} \int_{\mathbb{R}} h_n^\alpha(x,t)(\omega') \mathbf{P}'(d\omega')\, dx$$
$$= (\sigma C_H^{n/2})^\alpha \mathbf{E}' \int_{\mathbb{R}} h_n^\alpha(x,t)\, dx$$

$$\leq \int_{\mathbb{R}} \left[ \mathbf{E}' \left( \int_0^t \frac{p_{s^{2H}}(x/\sigma)}{s^{nH}} H_n\left(\frac{x/\sigma}{s^H}\right) I_n(K_H(s,\cdot)^{\otimes n}) \, ds \right)^2 \right]^{\alpha/2} dx$$

$$\leq \int_{\mathbb{R}} \left\{ \int_0^t \frac{p_{s^{2H}}(x/\sigma)}{s^{nH}} \left| H_n\left(\frac{x/\sigma}{s^H}\right) \right| [\mathbf{E}'(I_n(K_H(s,\cdot)^{\otimes n}))^2]^{1/2} \, ds \right\}^{\alpha} dx.$$

By (6.1)

$$\mathbf{E}'(I_n(K_H(s,\cdot)^{\otimes n}))^2 = n! \|K_H(s,\cdot)^{\otimes n}\|_2^2 = n! \|K_H(s,\cdot)\|_2^{2n} = n! C_H^n (s^{2H})^n.$$

Therefore,

$$\frac{\sigma^{\alpha}}{(n!)^{\alpha/2}} \mathbf{E}' \int_{\mathbb{R}} h_n^{\alpha}(x,t) \, dx \leq \int_{\mathbb{R}} \left( \int_0^t p_{s^{2H}}(x/\sigma) \left| H_n\left(\frac{x/\sigma}{s^H}\right) \right| ds \right)^{\alpha} dx$$

$$= \int_{\mathbb{R}} \left( \int_0^t p_1\left(\frac{x}{\sigma s^H}\right) \left| H_n\left(\frac{x/\sigma}{s^H}\right) \right| \frac{ds}{s^H} \right)^{\alpha} dx.$$

Observe that the function $\varphi(x) = p_1(x)^{1/2} |H_n(x)|$ is continuous and bounded on the entire real line. Therefore,

$$\int_0^t p_1\left(\frac{x}{\sigma s^H}\right) \left| H_n\left(\frac{x/\sigma}{s^H}\right) \right| \frac{ds}{s^H} \leq c \left( p_1\left(\frac{x}{\sigma t^H}\right) \right)^{1/2} \int_0^t \frac{ds}{s^H}$$

for some finite positive $c$ independent of $x$. Therefore,

(6.4) $$\mathbf{E}' \int_{\mathbb{R}} h_n^{\alpha}(x,t) \, dx < \infty,$$

and so each $(W_n(t), t \geq 0)$ is a well-defined S$\alpha$S process.

The next step is to check that each process $(W_n(t), t \geq 0)$ is self-similar, with the exponent of self-similarity given by (3.2). We will use two simple scaling facts. The first is simply

(6.5) $$K_H(au, w) = a^{H-1/2} K_H(u, w/a)$$

for all $a > 0$ and all $u, w$, and the second is a consequence of the self-similarity of a Brownian motion: for any $n \geq 1, m \geq 1$, any square-integrable symmetric functions $f_1, \ldots, f_m$ and $a > 0$,

(6.6) $$(I_n(f_i(\cdot a)), i = 1, \ldots, m) \stackrel{d}{=} (a^{-n/2} I_n(f_i), i = 1, \ldots, m).$$

We assume, for simplicity, that $\sigma = 1$.

Now let $m \geq 1$, $0 < t_1 < \cdots < t_m$, $\theta_1, \ldots, \theta_m \in \mathbb{R}$ and $a > 0$. We have

$$-\log \mathbf{E} \exp\left\{ iC_H^{n/2} \sum_{j=1}^m \theta_j W_n(at_j) \right\}$$

$$= \int_{\mathbb{R}} \left( \mathbf{E}' \left| \sum_{j=1}^m \theta_j \int_0^{at_j} \frac{p_{s^{2H}}(x)}{s^{nH}} H_n\left(\frac{x}{s^H}\right) I_n(K_H(s,\cdot)^{\otimes n}) \, ds \right|^{\alpha} \right) dx$$

$$= a^{\alpha-\alpha nH} \int_{\mathbb{R}} \left( \mathbf{E}' \left| \sum_{j=1}^{m} \theta_j \int_0^{t_j} \frac{p_{a^{2H}u^{2H}}(x)}{u^{nH}} H_n\left(\frac{x}{a^H u^H}\right) \right.\right.$$
$$\left.\left. \times I_n(K_H(au,\cdot)^{\otimes n})\,du \right|^\alpha \right) dx$$

$$= a^{\alpha-\alpha nH} \int_{\mathbb{R}} \left( \mathbf{E}' \left| \sum_{j=1}^{m} \theta_j \int_0^{t_j} \frac{p_{u^{2H}}(x/a^H)}{a^H u^{nH}} H_n\left(\frac{x}{a^H u^H}\right) \right.\right.$$
$$\left.\left. \times I_n(K_H(au,\cdot)^{\otimes n})\,du \right|^\alpha \right) dx$$

$$= a^{\alpha-\alpha nH+H-\alpha H} \int_{\mathbb{R}} \left( \mathbf{E}' \left| \sum_{j=1}^{m} \theta_j \int_0^{t_j} \frac{p_{u^{2H}}(y)}{u^{nH}} H_n\left(\frac{y}{u^H}\right) \right.\right.$$
$$\left.\left. \times I_n(K_H(au,\cdot)^{\otimes n})\,du \right|^\alpha \right) dy$$

$$= a^{\alpha-\alpha n/2+H-\alpha H} \int_{\mathbb{R}} \left( \mathbf{E}' \left| \sum_{j=1}^{m} \theta_j \int_0^{t_j} \frac{p_{u^{2H}}(y)}{u^{nH}} H_n\left(\frac{y}{u^H}\right) \right.\right.$$
$$\left.\left. \times I_n\left(K_H\left(u,\frac{\cdot}{a}\right)^{\otimes n}\right)du \right|^\alpha \right) dy$$

$$= a^{\alpha+H-\alpha H} \int_{\mathbb{R}} \left( \mathbf{E}' \left| \sum_{j=1}^{m} \theta_j \int_0^{t_j} \frac{p_{u^{2H}}(y)}{u^{nH}} H_n\left(\frac{y}{u^H}\right) I_n(K_H(u,\cdot)^{\otimes n})\,du \right|^\alpha \right) dy$$

$$= -\log \mathbf{E}\exp\left\{ iC_H^{n/2} \sum_{j=1}^{m} \theta_j a^{H'} W_n(at_j) \right\},$$

where the fifth equality follows from (6.5) and the sixth equality follows from (6.6). This proves the claimed self-similarity.

Because of self-similarity, it is enough to prove absolute continuity on the interval $[0,1]$. The proof will be done in three different cases; each consists of checking the conditions of Theorem 11.7.4 in [34].

If $0 < \alpha < 1$, we need to check that

$$\int_{\mathbb{R}} \mathbf{E}' \left( \int_0^1 \left| \frac{\partial h_n}{\partial t}(x,t) \right| dt \right)^\alpha dx < \infty.$$

This is, however, an immediate consequence of the computation leading to (6.4).

If $1 < \alpha < 2$, we need to check

(6.7) $$\int_0^1 \left( \int_{\mathbb{R}} \mathbf{E}' \left| \frac{\partial h_n}{\partial t}(x,t) \right|^\alpha dx \right)^{1/\alpha} dt < \infty.$$

We have, for $t > 0$

$$\int_{\mathbb{R}} \mathbf{E}'\left|\frac{\partial h_n}{\partial t}(x,t)\right|^\alpha dx = \int_{\mathbb{R}} \mathbf{E}'\left|\frac{p_{t^{2H}}(x)}{t^{nH}} H_n\left(\frac{x}{t^H}\right) I_n(K_H(t,\cdot)^{\otimes n})\right|^\alpha dx$$

$$\leq t^{-\alpha nH}(\mathbf{E}' I_n(K_H(t,\cdot)^{\otimes n})^2)^{\alpha/2}$$

$$\times \int_{\mathbb{R}} \frac{(p_1(x/t^H)|H_n(x/t^H)|)^\alpha}{t^{\alpha H}} dx$$

$$= b_n t^{-(\alpha-1)H},$$

for some $0 < b_n < \infty$. Since

$$\int_0^1 t^{-(\alpha-1)H/\alpha} dt < \infty,$$

(6.7) follows.

Finally, in the case $\alpha = 1$, the necessary and sufficient conditions for absolute continuity are less convenient to check. However, a stronger statement, that the process is absolutely continuous with a derivative in $L^p[0,1]$ for some $1 < p \leq 2$, requires checking that

$$\int_{\mathbb{R}} \mathbf{E}'\left(\int_0^1 \left|\frac{\partial h_n}{\partial t}(x,t)\right|^p dt\right)^{1/p} dx < \infty,$$

which follows in the same way as (6.4). We omit the repetitive details.

It remains to prove that the sequence (6.2) converges in probability. By Proposition 4 of [9], for every $x$ and $t$

$$l(x,t) = \sum_{n=0}^\infty h_n(x,t), \qquad \mathbf{P}'\text{-a.s.}$$

By definition (3.1) of the process $(Y(t), t \geq 0)$, it is enough to prove that

(6.8) $$\int_{\mathbb{R}} \mathbf{E}'\left|l(x,t) - \sum_{n=0}^m h_n(x,t)\right|^\alpha dx \to 0 \quad \text{as } m \to \infty.$$

We will estimate the expectation in (6.8) in two different ways. First note that, for every $m$,

$$E'\left|l(x,t) - \sum_{n=0}^m h_n(x,t)\right|^\alpha \leq C_\alpha\left[E' l(x,t)^\alpha + \left(E'\left(\sum_{n=0}^m h_n(x,t)\right)^2\right)^{\alpha/2}\right]$$

$$\leq C_\alpha[E' l(x,t)^\alpha + (E' l(x,t)^2)^{\alpha/2}]$$

$$\leq C_\alpha (E' l(x,t)^2)^{\alpha/2},$$

where $C_\alpha$ is a finite positive constant depending only on $\alpha$ and allowed to change from place to place. The argument used in (3.3) shows that

$$\int_{\mathbb{R}} (E'l(x,t)^2)^{\alpha/2}\,dx < \infty.$$

Therefore, given $\epsilon > 0$, one can choose $M \in (0, \infty)$ such that for all $m \geq 1$,

$$\int_{\mathbb{R}} \mathbf{E}'\left|l(x,t) - \sum_{n=0}^{m} h_n(x,t)\right|^\alpha dx \leq \epsilon + \int_{-M}^{M} \mathbf{E}'\left|l(x,t) - \sum_{n=0}^{m} h_n(x,t)\right|^\alpha dx.$$
(6.9)

Next, we estimate the expectation in (6.8) in a different way. Note that by the orthogonality of $(h_n)$'s with different $n$,

$$E'\left|l(x,t) - \sum_{n=0}^{m} h_n(x,t)\right|^\alpha \leq \left(E'\left(l(x,t) - \sum_{n=0}^{m} h_n(x,t)\right)^2\right)^{\alpha/2}$$
$$= \left(E'\left(\sum_{n=m+1}^{\infty} h_n(x,t)\right)^2\right)^{\alpha/2}$$
$$= \left(\sum_{n=m+1}^{\infty} E'(h_n(x,t))^2\right)^{\alpha/2}$$

and, as in the proof of Proposition 4 in [9], we conclude that there exist some $\delta_{m,t} \to 0$ as $m \to \infty$ such that, for all $x \in \mathbb{R}$,

$$E'\left|l(x,t) - \sum_{n=0}^{m} h_n(x,t)\right|^\alpha \leq \delta_{m,t}.$$

Substituting this bound into (6.9), we conclude that

$$\int_{\mathbb{R}} \mathbf{E}'\left|l(x,t) - \sum_{n=0}^{m} h_n(x,t)\right|^\alpha dx \leq \epsilon + 2M\delta_{m,t}.$$

Letting first $m \to \infty$ and then $\epsilon \to 0$ proves (6.8), and so the proof of the theorem is now complete. □

REMARK 6.2. It is clear from the proof of the theorem that the derivative of each process $(W_n(t), t \geq 0)$ is in $L^p[0,1]$ for a range of $p > 1$ in all cases, and not only for $\alpha = 1$. We will not pursue this point here, however.

**7. Convergence of the random reward scheme.** In this section we establish the limit theorem in the random reward scheme discussed in the Introduction. We start by setting up the notation. Let $(W_k^{(i)}, k \in \mathbf{Z}, i \geq 1)$ be an array of i.i.d. symmetric random variables whose distribution satisfies (1.1). Let $(V_k^{(i)}, k \geq 1, i \geq 1)$ be an array of i.i.d. mean zero and unit variance integer-valued random variables, independent of $(W_k^{(i)}, k \in \mathbf{Z}, i \geq 1)$. Let $S_n^{(i)} = V_1^{(i)} + \cdots + V_n^{(i)}$, $n \geq 0$, be the $i$th random walk, $i = 1, 2, \ldots$, and define for $j \in \mathbf{Z}$ and $n \geq 1$

$$\varphi(j, n; i) = \sum_{k=1}^{n} \mathbf{1}(S_k^{(i)} = j) \tag{7.1}$$

to be the number of times the $i$th random walk visits the state $j$ by time $n$, $i = 1, 2, \ldots$. Define $\varphi(j, t; i)$ for noninteger values of $t \geq 0$ by interpolating linearly between $\varphi(j, n; i)$ and $\varphi(j, n+1; i)$ if $n \leq t < n+1$ [we use $\varphi(j, 0; i) = 0$]. Notice that the total reward earned by the $i$th user by time $t$ can be written as

$$R^{(i)}(t) = \sum_{k=-\infty}^{\infty} W_k^{(i)} \varphi(k, t; i)$$

(of course, this is really a linear interpolation for noninteger $t$). The limit theorem below shows that, if both the number of users and the time scale grow *at an arbitrary rate*, then the properly normalized total reward converges weakly to the FBM-1/2-local time fractional symmetric $\alpha$-stable motion. This is related to the convergence result in [15] (which allows more general random walks) where only one user is present.

THEOREM 7.1. *For every sequence $(b_n)$ of positive integers with $b_n \to \infty$ we have, as $n \to \infty$,*

$$\left( \frac{1}{(nb_n^{(\alpha+1)/2})^{1/\alpha}} \sum_{i=1}^{n} R^{(i)}(b_n t), t \geq 0 \right) \implies ((2/C_\alpha)^{1/\alpha} \sigma_W Y(t), t \geq 0) \tag{7.2}$$

*weakly in $C([0, \infty))$, where $(Y(t), t \geq 0)$ is the FBM-1/2-local time fractional symmetric $\alpha$-stable motion defined in* (3.1) *(with the local time being that of a standard Brownian motion). Here, $\sigma_W$ is the tail weight in* (1.1) *and $C_\alpha$ is the stable tail constant given by*

$$C_\alpha = \left( \int_0^\infty x^{-\alpha} \sin x \, dx \right)^{-1}. \tag{7.3}$$

PROOF. By extending the probability space on which random objects are defined, if necessary, we can construct a sequence of i.i.d. standard Brownian motions $(B^{(i)}(t), t \geq 0)$, $i = 1, 2, \ldots$, with jointly continuous local time processes $(l^{(i)}(x, t), t \geq 0, t \in \mathbb{R})$, $i = 1, 2, \ldots$, such that, for every $T > 0$,

$$\sup_{x \in \mathbb{R}, 0 \leq t \leq nT} \left| \varphi([x], t; i) - n^{1/2} l^{(i)}\left( \frac{x}{\sqrt{n}}, \frac{t}{n} \right) \right| \to 0 \tag{7.4}$$

in probability as $n \to \infty$, $i = 1, 2, \ldots$; see Theorem 1 of [5]. Define, for $n \geq 1$,

$$(7.5) \quad X_n(t) = \frac{1}{(nb_n^{1/2})^{1/\alpha}} \sum_{i=1}^{n} \sum_{k=-\infty}^{\infty} W_k^{(i)} l^{(i)}\left(\frac{k}{\sqrt{b_n}}, t\right), \quad t \geq 0.$$

Notice that, for $t \geq 0$,

$$E_n(t) := \frac{1}{(nb_n^{(\alpha+1)/2})^{1/\alpha}} \sum_{i=1}^{n} R^{(i)}(b_n t) - X_n(t)$$

$$(7.6) \quad = \frac{1}{(nb_n^{(\alpha+1)/2})^{1/\alpha}}$$

$$\times \sum_{i=1}^{n} \sum_{k=-\infty}^{\infty} W_k^{(i)}\left(\varphi(k, b_n t; i) - b_n^{1/2} l^{(i)}\left(\frac{k}{\sqrt{b_n}}, t\right)\right).$$

We first prove that, for every $t > 0$,

$$(7.7) \quad E_n(t) \to 0 \quad \text{in probability.}$$

For notational simplicity, we prove (7.7) for $t = 1$.

First, it follows from the tail behavior (1.1) that there is a constant $b > 0$ such that

$$(7.8) \quad |W_k^{(i)}| \stackrel{\text{st}}{\leq} b(1 + |R_k^{(i)}|)$$

(in the sense of stochastic comparison), where $(R_k^{(i)}, k \in \mathbf{Z}, i \geq 1)$ is an array of i.i.d. standard S$\alpha$S random variables. Therefore, by the contraction inequality (see Section 1.2 of [18]) we conclude that

$$P(|E_n(1)| > \epsilon)$$

$$\leq 2P\Bigg(\frac{1}{(nb_n^{(\alpha+1)/2})^{1/\alpha}}$$

$$\times \sum_{i=1}^{n} \sum_{k=-\infty}^{\infty} \varepsilon_k^{(i)}(1 + |R_k^{(i)}|)$$

$$(7.9) \quad \times \left(\varphi(k, b_n; i) - b_n^{1/2} l^{(i)}\left(\frac{k}{\sqrt{b_n}}, 1\right)\right) > \epsilon/b\Bigg)$$

$$\leq 2P\Bigg(\frac{1}{(nb_n^{(\alpha+1)/2})^{1/\alpha}}$$

$$\times \sum_{i=1}^{n} \sum_{k=-\infty}^{\infty} \varepsilon_k^{(i)}\left(\varphi(k, b_n; i) - b_n^{1/2} l^{(i)}\left(\frac{k}{\sqrt{b_n}}, 1\right)\right) > \epsilon/(2b)\Bigg)$$

$$+ 2P\left(\frac{1}{(nb_n^{(\alpha+1)/2})^{1/\alpha}}\right.$$
$$\left.\times \sum_{i=1}^{n}\sum_{k=-\infty}^{\infty} R_k^{(i)}\left(\varphi(k,b_n;i) - b_n^{1/2}l^{(i)}\left(\frac{k}{\sqrt{b_n}},1\right)\right) > \epsilon/(2b)\right)$$
$$:= p_1(n) + p_2(n),$$

where $(\varepsilon_k^{(i)}, k \in \mathbf{Z}, i \geq 1)$ is an array of i.i.d standard symmetric Rademacher random variables. We need to show that

(7.10) $\qquad p_j(n) \to 0 \quad \text{as } n \to \infty \text{ for } j = 1, 2.$

We estimate $p_2(n)$. Note that

$$p_2(n)/2$$
$$= P\left(\left(\frac{1}{nb_n^{(\alpha+1)/2}}\right.\right.$$
$$\left.\left.\times \sum_{i=1}^{n}\sum_{k=-\infty}^{\infty}\left|\varphi(k,b_n;i) - b_n^{1/2}l^{(i)}\left(\frac{k}{\sqrt{b_n}},1\right)\right|^\alpha\right)^{1/\alpha} R_1^{(1)} > \epsilon/(2b)\right),$$

and so the statement (7.10) with $j = 2$ will follow once we show that

(7.11) $\qquad \dfrac{1}{nb_n^{(\alpha+1)/2}} \sum_{i=1}^{n}\sum_{k=-\infty}^{\infty}\left|\varphi(k,b_n;i) - b_n^{1/2}l^{(i)}\left(\dfrac{k}{\sqrt{b_n}},1\right)\right|^\alpha \to 0$

in probability as $n \to \infty$. The expectation of the expression in the left-hand side of (7.11) is

$$\frac{1}{b_n^{(\alpha+1)/2}} E \sum_{k=-\infty}^{\infty}\left|\varphi(k,b_n;1) - b_n^{1/2}l^{(1)}\left(\frac{k}{\sqrt{b_n}},1\right)\right|^\alpha$$

$$\leq \frac{1}{b_n^{(\alpha+1)/2}} E \sum_{k=-\infty}^{\infty}\left|\varphi(k,b_n;1) - b_n^{1/2}l^{(1)}\left(\frac{k}{\sqrt{b_n}},1\right)\right|^\alpha$$
$$\times \mathbf{1}\left(\left|\varphi(k,b_n;1) - b_n^{1/2}l^{(1)}\left(\frac{k}{\sqrt{b_n}},1\right)\right| \leq 1\right)$$
$$+ \frac{1}{b_n^{(\alpha+1)/2}}$$
$$\times E\left[\left(\sum_{k=-\infty}^{\infty}\left|\varphi(k,b_n;1) - b_n^{1/2}l^{(1)}\left(\frac{k}{\sqrt{b_n}},1\right)\right|^2\right)^{\alpha/2}\right]$$

$$\times \left( \sum_{k=-\infty}^{\infty} \mathbf{1}\left( \left|\varphi(k,b_n;1) - b_n^{1/2} l^{(1)}\left(\frac{k}{\sqrt{b_n}},1\right)\right| > 1 \right) \right)^{1-\alpha/2} \right]$$

$$\leq \frac{1}{b_n^{(\alpha+1)/2}} E \sum_{k=-\infty}^{\infty} \left|\varphi(k,b_n;1) - b_n^{1/2} l^{(1)}\left(\frac{k}{\sqrt{b_n}},1\right)\right|^{\alpha}$$

$$\times \mathbf{1}\left( \left|\varphi(k,b_n;1) - b_n^{1/2} l^{(1)}\left(\frac{k}{\sqrt{b_n}},1\right)\right| \leq 1 \right)$$

$$+ \frac{1}{b_n^{(\alpha+1)/2}} E \left[ \left( \sum_{k=-\infty}^{\infty} \left|\varphi(k,b_n;1) - b_n^{1/2} l^{(1)}\left(\frac{k}{\sqrt{b_n}},1\right)\right|^2 \right)^{\alpha/2} \right.$$

$$\times (2M^{(1)}(b_n)+1)^{1-\alpha/2}$$

$$\left. \times \mathbf{1}\left( \sup_{k\in\mathbf{Z}} \left|\varphi(k,b_n;1) - b_n^{1/2} l^{(1)}\left(\frac{k}{\sqrt{b_n}},1\right)\right| > 1 \right) \right]$$

$$:= p_{21}(n) + p_{22}(n),$$

where

$$M^{(i)}(m) = \max\left( \sup_{0\leq k\leq m} |S_k^{(i)}|, \sqrt{m} \sup_{0\leq t\leq 1} |B^{(i)}(t)| \right).$$

For the second inequality above, we have bounded a sum from above by the number of nonvanishing terms times the largest nonvanishing term. A similar argument will be used in the sequel without further comment. We have

$$p_{21}(n) \leq \frac{1}{b_n^{(\alpha+1)/2}} E(2M^{(1)}(b_n)+1)$$

$$\leq c \frac{1}{b_n^{(\alpha+1)/2}} b_n^{1/2} = c b_n^{-\alpha/2} \to 0 \qquad \text{as } n \to \infty.$$

Furthermore,

$$p_{22}(n) \leq \frac{1}{b_n^{(\alpha+1)/2}} \left( E \sum_{k=-\infty}^{\infty} \left|\varphi(k,b_n;1) - b_n^{1/2} l^{(1)}\left(\frac{k}{\sqrt{b_n}},1\right)\right|^2 \right)^{\alpha/2}$$

$$\times (E(2M^{(1)}(b_n)+1)\mathbf{1}(\Delta_1(b_n)>1))^{1-\alpha/2},$$

where

$$\Delta_i(n) = \sup_{k\in\mathbf{Z}} \left|\varphi(k,b_n;i) - b_n^{1/2} l^{(i)}\left(\frac{k}{\sqrt{b_n}},1\right)\right|.$$

Using Lemma 1 of [15] and the fact that the largest value of a Brownian local time at time 1 has all moments finite, the first expectation in the right-hand side is bounded above by $cb_n^{3/2}$. Therefore,

$$p_{22}(n) \leq c\frac{1}{b_n^{(\alpha+1)/2}}b_n^{3\alpha/4}(EM^{(1)}(b_n)^{3/2})^{(2-\alpha)/3}(P(\Delta_1(b_n) > 1))^{(1-\alpha/2)/3}$$

$$\leq c\frac{1}{b_n^{(\alpha+1)/2}}b_n^{3\alpha/4}(b_n^{3/4})^{(2-\alpha)/3}(P(\Delta_1(b_n) > 1))^{(1-\alpha/2)/3}$$

$$\leq c(P(\Delta_1(b_n) > 1))^{(1-\alpha/2)/3} \to 0 \quad \text{as } n \to \infty,$$

by (7.4) (as always, $c$ is a finite positive constant that may change from instance to instance). Therefore, (7.11) holds, and so we have established (7.10) for $j = 2$. The proof for $j = 1$ is similar. Thus, we have obtained (7.7).

The next step is to show that the finite-dimensional distributions of the process $(X_n(t), t \geq 0)$ in (7.5) converge to those of $(Y(t), t \geq 0)$. For this, it is enough to show that, for every $k \geq 1$, $0 < t_1 < \cdots < t_k$ and $\theta_1, \ldots, \theta_k \in \mathbb{R}$,

$$\sum_{j=1}^{k} \theta_j X_n(t_j) \implies \sum_{j=1}^{k} \theta_j Y(t_j) \quad \text{as } n \to \infty.$$

We will see that this is true for $k = 1$ and $t_1 = 1$; the general case is only notationally different. That is, we will show that

(7.12) $$\frac{1}{(nb_n^{1/2})^{1/\alpha}} \sum_{i=1}^{n} \sum_{k=-\infty}^{\infty} W_k^{(i)} l^{(i)}\left(\frac{k}{\sqrt{b_n}}, 1\right) \implies Y(1) \quad \text{as } n \to \infty.$$

By Theorem 8 in Chapter 6 of [24] it is enough to prove that, for every $\lambda > 0$,

(7.13)
$$\lim_{n \to \infty} nP\left(\sum_{k=-\infty}^{\infty} W_k^{(1)} l^{(1)}\left(\frac{k}{\sqrt{b_n}}, 1\right) > \lambda(nb_n^{1/2})^{1/\alpha}\right)$$
$$\to \sigma_W^\alpha E \int_{\mathbb{R}} l(x,t)^\alpha \, dx \, \lambda^{-\alpha}$$

and

(7.14)
$$\lim_{\epsilon \to 0} \limsup_{n \to \infty} \frac{n}{(nb_n^{1/2})^{2/\alpha}}$$
$$\times E\left[\left(\sum_{k=-\infty}^{\infty} W_k^{(1)} l^{(1)}\left(\frac{k}{\sqrt{b_n}}, 1\right)\right)^2\right.$$
$$\left.\times \mathbf{1}\left(\left|\sum_{k=-\infty}^{\infty} W_k^{(1)} l^{(1)}\left(\frac{k}{\sqrt{b_n}}, 1\right)\right| \leq \epsilon(nb_n^{1/2})^{1/\alpha}\right)\right] = 0$$

(we have used the symmetry of $W$'s to simplify the conditions).

We start by checking (7.13). The first step is to prove that, for every $\lambda > 0$,

$$(7.15) \quad \lim_{K \to \infty} \limsup_{n \to \infty} nP\left(\sum_{|k|>K\sqrt{b_n}} W_k^{(1)} l^{(1)}\left(\frac{k}{\sqrt{b_n}}, 1\right) > \lambda(nb_n^{1/2})^{1/\alpha}\right) = 0.$$

By using the contraction inequality, the stochastic comparison (7.8) and the notation following it, it is enough to prove that, for every $\lambda > 0$, (7.15) holds with each $W_k^{(1)}$ being replaced by $R_k^{(1)}$, and with each $W_k^{(1)}$ being replaced by $\varepsilon_k^{(1)}$. The two statements are similar; we only present the argument in the case of stable weights. In that case, the expression corresponding to that in the left-hand side of (7.15) is equal to

$$nP\left(\left(\sum_{|k|>K\sqrt{b_n}} l^{(1)}\left(\frac{k}{\sqrt{b_n}}, 1\right)^\alpha\right)^{1/\alpha} R_1^{(1)} > \lambda(nb_n^{1/2})^{1/\alpha}\right)$$

and, for some positive constant $c$, this bounded from above by

$$n\left(c\lambda^{-\alpha}(nb_n^{1/2})^{-1} E\left[\sum_{|k|>K\sqrt{b_n}} l^{(1)}\left(\frac{k}{\sqrt{b_n}}, 1\right)^\alpha\right]\right)$$

$$= c\lambda^{-\alpha} b_n^{-1/2} E\left[\sum_{|k|>K\sqrt{b_n}} l^{(1)}\left(\frac{k}{\sqrt{b_n}}, 1\right)^\alpha\right]$$

$$\leq c\lambda^{-\alpha} E\left(\sup_{x \in \mathbb{R}} l^{(1)}(x,1)\right)^\alpha \sum_{|k|>K\sqrt{b_n}} P\left(\sup_{0 \leq s \leq 1} |B^{(1)}(s)| \geq \frac{k}{\sqrt{b_n}}\right)$$

$$\to 2c\lambda^{-\alpha} \int_K^\infty P\left(\sup_{0 \leq s \leq 1} |B^{(1)}(s)| > x\right) dx.$$

Since the final expression converges to 0 as $K \to \infty$, we have (7.15).

Now fix $K$ and $\lambda > 0$. The usual "largest jump" large deviations approach (see, e.g., [22]) and the continuity of the local time give us that, as $n \to \infty$,

$$nP\left(\sum_{|k| \leq K\sqrt{b_n}} W_k^{(1)} l^{(1)}\left(\frac{k}{\sqrt{b_n}}, 1\right) > \lambda(nb_n^{1/2})^{1/\alpha}\right)$$

$$\sim nP\left(\max_{|k| \leq K\sqrt{b_n}} W_k^{(1)} l^{(1)}\left(\frac{k}{\sqrt{b_n}}, 1\right) > \lambda(nb_n^{1/2})^{1/\alpha}\right)$$

$$(7.16) \quad \sim n \sum_{|k| \leq K\sqrt{b_n}} P\left(W_k^{(1)} l^{(1)}\left(\frac{k}{\sqrt{b_n}}, 1\right) > \lambda(nb_n^{1/2})^{1/\alpha}\right)$$

$$\sim n \int_{-K\sqrt{b_n}}^{K\sqrt{b_n}} P\left(W_k^{(1)} l^{(1)}\left(\frac{x}{\sqrt{b_n}}, 1\right) > \lambda (nb_n^{1/2})^{1/\alpha}\right) dx$$

$$= \int_{-K}^{K} (nb_n^{1/2}) P(W_k^{(1)} l^{(1)}(y, 1) > \lambda (nb_n^{1/2})^{1/\alpha}) \, dy$$

$$\to \int_{-K}^{K} \sigma_W^\alpha \lambda^{-\alpha} E(l^{(1)}(y, 1))^\alpha \, dy$$

(see, e.g., (2.7) in [31]).

(7.13) now follows from (7.15) and (7.16).

To show (7.14), note that

$$\frac{n}{(nb_n^{1/2})^{2/\alpha}} E\left[\left(\sum_{k=-\infty}^{\infty} W_k^{(1)} l^{(1)}\left(\frac{k}{\sqrt{b_n}}, 1\right)\right)^2\right.$$

$$\left. \times \mathbf{1}\left(\left|\sum_{k=-\infty}^{\infty} W_k^{(1)} l^{(1)}\left(\frac{k}{\sqrt{b_n}}, 1\right)\right| \leq \epsilon (nb_n^{1/2})^{1/\alpha}\right)\right]$$

$$\leq n \int_0^{\epsilon^2} P\left[\left|\sum_{k=-\infty}^{\infty} W_k^{(1)} l^{(1)}\left(\frac{k}{\sqrt{b_n}}, 1\right)\right| > x^{1/2} b_n^{1/2\alpha} n^{1/\alpha}\right] dx.$$

Using stochastic domination and contraction principle as above allows us to replace the random variables $W_k^{(i)}$ in the above expression by S$\alpha$S random variables and by Rademacher random variables and, as before, we only consider the former (because they have heavier tails). In that case, we have

$$n \int_0^{\epsilon^2} P\left[\left|\sum_{k=-\infty}^{\infty} R_k^{(1)} l^{(1)}\left(\frac{k}{\sqrt{b_n}}, 1\right)\right| > x^{1/2} b_n^{1/2\alpha} n^{1/\alpha}\right] dx$$

$$= n \int_0^{\epsilon^2} P\left[|R_1^{(1)}|\left(\sum_{k=-\infty}^{\infty} \left(l^{(1)}\left(\frac{k}{\sqrt{b_n}}, 1\right)\right)^\alpha\right)^{1/\alpha} > x^{1/2} b_n^{1/2\alpha} n^{1/\alpha}\right] dx$$

$$\leq c b_n^{-1/2} E \sum_{k=-\infty}^{\infty} \left(l^{(1)}\left(\frac{k}{\sqrt{b_n}}, 1\right)\right)^\alpha \int_0^{\epsilon^2} x^{-\alpha/2} \, dx,$$

from which (7.14) would follow once we check uniform boundedness of the $n$-dependent coefficient above. However, this follows from

$$b_n^{-1/2} E \sum_{k=-\infty}^{\infty} \left(l^{(1)}\left(\frac{k}{\sqrt{b_n}}, 1\right)\right)^\alpha \leq E\left[\sup_{x \in \mathbb{R}} (l(x,1))^\alpha \left(2 \sup_{0 \leq t \leq 1} |B^{(1)}(t)| + 1\right)\right] < \infty.$$

Therefore, we have (7.14) and, thus, convergence of the finite-dimensional distributions in (7.2).

It remains to prove tightness. Write, for $M > 0$,

$$\frac{1}{(nb_n^{(\alpha+1)/2})^{1/\alpha}} \sum_{i=1}^n R^{(i)}(b_n t)$$

(7.17)
$$= \frac{1}{(nb_n^{(\alpha+1)/2})^{1/\alpha}} \sum_{i=1}^n \sum_{k=-\infty}^{\infty} W_k^{(i)} \mathbf{1}(|W_k^{(i)}| > Mn^{1/\alpha} b_n^{1/2\alpha}) \varphi(k, b_n t; i)$$

$$+ \frac{1}{(nb_n^{(\alpha+1)/2})^{1/\alpha}} \sum_{i=1}^n \sum_{k=-\infty}^{\infty} W_k^{(i)} \mathbf{1}(|W_k^{(i)}| \leq Mn^{1/\alpha} b_n^{1/2\alpha}) \varphi(k, b_n t; i)$$

$$:= Y_n(t) + Z_n(t), \quad t \geq 0.$$

Notice that

$$P\left(\sup_{0 \leq t \leq 1} |Y_n(t)| > 0\right) \leq 1 - [P(\text{for all } |k| \leq M^{(1)}(b_n), |W_k^{(1)}| \leq Mn^{1/\alpha} b_n^{1/2\alpha})]^n.$$

Since for large $n$, with a changing constant $c$,

$$P(\text{for all } |k| \leq M^{(1)}(b_n), |W_k^{(1)}| \leq Mn^{1/\alpha} b_n^{1/2\alpha})$$
$$= E[P(|W_1^{(1)}| \leq Mn^{1/\alpha} b_n^{1/2\alpha})]^{2M^{(1)}(b_n)+1}$$
$$\geq E[1 - cM^{-\alpha} n^{-1} b_n^{-1/2}]^{2M^{(1)}(b_n)+1}$$
$$\geq E \exp\{-cM^{-\alpha} n^{-1} b_n^{-1/2} (2M^{(1)}(b_n) + 1)\},$$

we obtain, using the inequality $e^{-x} \geq 1 - x$ for $x \geq 0$ and maximal inequality for martingales,

$$P\left(\sup_{0 \leq t \leq 1} |Y_n(t)| > 0\right) \leq 1 - [1 - cM^{-\alpha} n^{-1} b_n^{-1/2} E(2M^{(1)}(b_n) + 1)]^n$$
$$\leq 1 - [1 - cM^{-\alpha} n^{-1}]^n \to 1 - \exp\{-cM^{-\alpha}\},$$

as $n \to \infty$.

Since the last expression converges to zero as $M \to \infty$, it follows from (7.17) that it is enough to prove that, for each fixed $M$, the process $(Z_n(t), 0 \leq t \leq 1)$ is tight.

However, for all $0 \leq s < t \leq 1$, we have

$$E(Z_n(t) - Z_n(s))^2 = \frac{1}{n^{2/\alpha - 1} b_n^{(\alpha+1)/\alpha}} E[(W_1^{(1)})^2 \mathbf{1}(|W_1^{(1)}| \leq Mn^{1/\alpha} b_n^{1/2\alpha})]$$
$$\times E \sum_{k=-\infty}^{\infty} (\varphi(k, b_n t; 1) - \varphi(k, b_n s; 1))^2.$$

Since, for large $x$,

$$E[(W_1^{(1)})^2 \mathbf{1}(|W_1^{(1)}| \leq x)] \leq 4 \int_0^{x^2} y P(W_1^{(1)} > y)\, dy \leq cx^{2-\alpha},$$

we see that, for large $n$,

$$E(Z_n(t) - Z_n(s))^2 \leq c b_n^{-3/2} E \sum_{k=-\infty}^{\infty} (\varphi(k, b_n t; 1) - \varphi(k, b_n s; 1))^2$$

$$\leq c(t-s)^{3/2}$$

as in the proof of Lemma 7 in [15]. We can now appeal to Theorem 12.3 in [4] to prove tightness of the family of the processes $(Z_n(t), 0 \leq t \leq 1)$ and, hence, complete the proof. □

**8. Discussion and possible extensions.** We briefly mention several issues related to the model constructed in this paper.

It is clear that self-similar S$\alpha$S processes with stationary increments could be constructed using local times of self-similar processes with stationary increments other than fractional Brownian motions. Symmetric stable Lévy motions with index of stability between 1 and 2 are obvious examples. One could also consider additive functionals other than local times.

For the random reward scheme considered in Section 7, it is clear that, in order to obtain in limit FBM-$H$-local time fractional symmetric $\alpha$-stable motion with $H \neq 1/2$, one has to introduce sufficiently long memory in the sequence of steps of each random walk $(V_k^{(i)}, k \geq 1)$. One way to do it is to take a stationary integer-valued sequence with slowly decaying correlations; alternatively, a certain reinforcement mechanism could be used. This is left for a future work.

It is also instructive to note that, in Section 7, one obtains the same limit regardless of how fast the number of users grows. However, if one considers instead (as is common in the literature) a *fluid input system*, where the random reward is not gained instantaneously, but, instead, obtained over a stretch of time, it is likely that different limits would be obtained, depending on the number of users. Possible limits there would, probably, include fractional Brownian motions, FBM-$H$-local time fractional symmetric $\alpha$-stable motions and, perhaps, additional limit processes.

**Acknowledgments.** The authors would like to thank Yimin Xiao for a useful discussion of the properties of local times and, in particular, of (2.9). A large part of this research was done during Samorodnitsky's visits to Université Paul Sabatier, whose support is gratefully acknowledged.

LABORATOIRE DE STATISTIQUE
ET DE PROBABILITÉS
UNIVERSITÉ PAUL SABATIER
INSTITUT DE MATHÉMATIQUE
118 ROUTE DE NARBONNE
31062 TOULOUSE
FRANCE
E-MAIL: scohen@cict.fr

SCHOOL OF OPERATIONS RESEARCH
AND INDUSTRIAL ENGINEERING
CORNELL UNIVERSITY
ITHACA, NEW YORK 14853
USA
E-MAIL: gennady@orie.cornell.edu